\documentclass[12pt]{article}
\usepackage[utf8]{inputenc}
\usepackage{amsmath}
\usepackage{amsfonts}
\usepackage{amssymb}
\usepackage{amsthm}
\usepackage{relsize}
\usepackage{array}
\usepackage{multirow}
\usepackage{tabu}
\usepackage{soul}
\usepackage{graphicx}
\usepackage{epigraph}
\usepackage{float}
\usepackage[english]{babel}
\usepackage[usenames, dvipsnames]{color}
\usepackage{fancyhdr}
\usepackage[Sonny]{fncychap}
\usepackage{tikz-cd}
\usepackage{tkz-euclide}
\usepackage{sseq}
\usepackage{newfloat}
\usepackage{amsmath,mathtools}
\usepackage{pgfplots}
\usepackage{geometry}
 \usepackage{sectsty}
\sectionfont{\fontsize{14}{15}\selectfont}

\usepackage[all]{xy}
\DeclareFloatingEnvironment[fileext=dia]{diagram}
\theoremstyle{definition}

\newtheorem{remark}{Remark}

\usepackage[most]{tcolorbox}

\makeatletter
\newtcolorbox{note}[1][]{%
	breakable,
	enhanced jigsaw, % better frame drawing
	borderline west={3pt}{0pt}{black!10!white}, % straigh vertical line at the left edge
	borderline south={1pt}{0pt}{black!10!white}, 
	borderline east={1pt}{0pt}{black!10!white},
	borderline north={1pt}{0pt}{black!10!white},
	sharp corners, % No rounded corners
	boxrule=0pt, % no real frame,
	attach title to upper, % Move the title into the box,
	left=0pt,
	right=0pt,
	top=0pt,
	bottom=0pt,
	boxsep=5pt,
	colback=white,
	frame hidden,
	#1
}
\makeatother

\makeatletter
\newtcolorbox{note1}[1][]{%
	breakable,
	enhanced jigsaw, % better frame drawing
	sharp corners, % No rounded corners
	boxrule=0pt, % no real frame,
	attach title to upper, % Move the title into the box,
	fontupper=\linespread{1.1}\fontfamily{qpl}\selectfont,
	fontlower=\linespread{1.1}\fontfamily{qpl}\selectfont, 
	left=0pt,
	right=0pt,
	top=0pt,
	bottom=0pt,
	boxsep=3pt,
	colback=green!3!white,
	frame hidden,
	before skip=10pt plus 2pt,after skip=10pt plus 2pt,
	#1
}
\makeatother

\makeatletter
\newcommand\tabfill[1]{%
	\dimen@\linewidth
	\advance\dimen@\@totalleftmargin
	\advance\dimen@-\dimen\@curtab
	\parbox[t]\dimen@{#1\ifhmode\strut\fi}%
}
\makeatother

\usepackage{tabularx}
\usepackage{imakeidx}
\usepackage{hyperref}
\hypersetup{colorlinks={true},linkcolor={blue},citecolor={blue},urlcolor={blue}, linktoc=all}  
 
 \usepackage[capitalize, nameinlink]{cleveref}
 \crefdefaultlabelformat{#2\textbf{#1}#3}
 \crefname{figure}{Figure}{Figures} 
 \Crefname{figure}{Figure}{Figures}
 \crefname{table}{Table}{Tables}
 \Crefname{table}{Table}{Tables}
 \crefname{section}{\S\hspace{-1mm}}{\S\hspace{-1mm}}
 \Crefname{section}{\S\hspace{-1mm}}{\S\hspace{-1mm}}
 \crefname{equation}{}{}
 \Crefname{equation}{}{}
 \crefname{example}{Geometric Pattern}{Geometric Patterns} 
 \Crefname{example}{Geometric Pattern}{Geometric Patterns}

 \usepackage{graphicx,tipa}% http://ctan.org/pkg/{graphicx,tipa}

 \usepackage{footnote}

\begin{document}

\title{\textbf{Similarity of Triangles and Intercept Theorem in Elamite Mathematics}}

\author{Nasser Heydari\footnote{Email: nasser.heydari@mun.ca}~ and  Kazuo Muroi\footnote{Email: edubakazuo@ac.auone-net.jp}}

\maketitle

\begin{abstract}
In this article, we study similarity of triangles in the Susa Mathematical Texts (\textbf{SMT}). We also suggest  that the Susa scribes were aware of intercept theory because they  used this theorem in solving a complicated system of equations. 
\end{abstract}

\section{Introduction}
Applications of similarity of triangles occur in some texts of the \textbf{SMT} such as \textbf{SMT No.\,18}, \textbf{SMT No.\,23}, and \textbf{SMT No.\,25}. We have  already mentioned  this technique in   \textbf{SMT No.\,23} and \textbf{SMT No.\,25} whose subject is the bisection of a trapezoid by a transversal line. For a full discussion about the mathematical interpretations of these two texts, see \cite{HM22-2}.

Here, we carefully examine   \textbf{SMT No.\,18}, which solves a complicated system of equations, and give our mathematical interpretation. This text  was inscribed by an Elamite scribe  between 1894--1595 BC on one of  26 clay tablets excavated from Susa in  southwest Iran by French archaeologists in 1933. The texts of all the tablets,  along with their interpretations, were first published in 1961 (see \cite{BR61}).

\textbf{SMT No.\,18}\footnote{The reader can see   this tablet on the website of the Louvre's collection. Please see \url{https://collections.louvre.fr/en/ark:/53355/cl010186428} for obverse  and   reverse.} contains only a single problem which shows one of the characteristics of Babylonian mathematics known as indifference to dimensions\footnote{Babylonians added  lengths to   areas or   volumes, which have different units. It seems that they were unconcerned by this practice. On the other hand,  the Greeks only added lengths to lengths or areas to areas and so on.}. In fact, the product of an area multiplied by another area, which would be meaningless to the ancient Greeks, occurs in one of the three equations given in this text. Moreover,   the scribe of this tablet has great skill,    introducing new variables in order to solve a system of complex simultaneous equations.

\section{Similarity   of Triangles} 
\begin{figure}[H]
	\centering
	\includegraphics[scale=1]{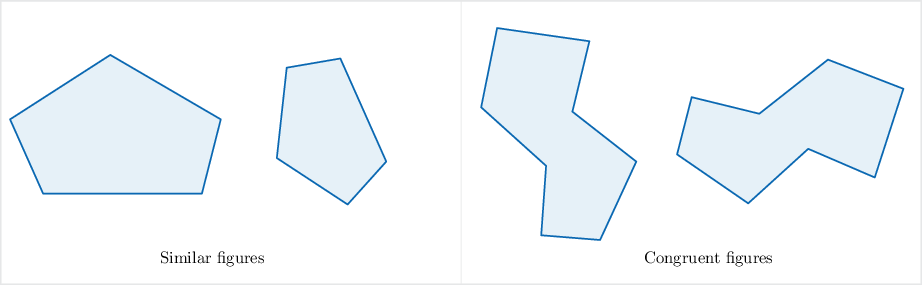}
	%\vspace{-1mm}
	\caption{Similar  and congruent figures}
	\label{Figure1-1}
\end{figure}

Two figures in the plane are said to be \textit{similar} if they can be obtained from one another  by applying    a combination of the following transformations: scaling, translating, rotating or reflecting. For example, all circles and regular $n$-gons for $n\geq 3$  are similar to each other, but this is not true for general isosceles triangles. We usually use the symbol ``$\sim$'' for similarity. In particular, if two figures are similar and the scaling coefficient is 1, they are called \textit{congruent} (see \cref{Figure1-1}).

\begin{figure}[H]
	\centering
	\includegraphics[scale=1]{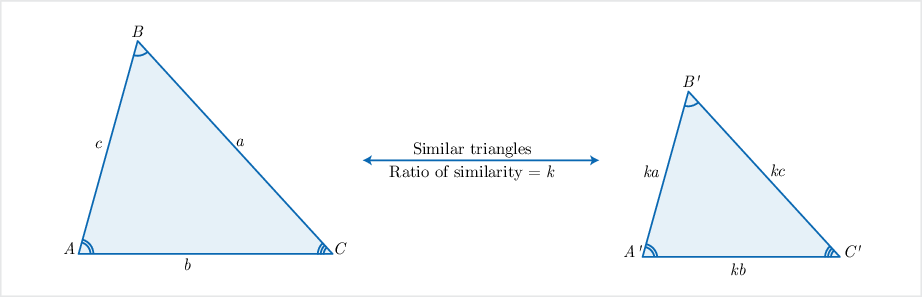}
	%\vspace{-1mm}
	\caption{Similar  triangles}
	\label{Figure1}
\end{figure}

Among all  polygons, the similarity of triangles has been of great interest to mathematicians, which has led to  different theorems giving the conditions under which two triangles are similar. Since  similarity definitions are equivalent,  no definition takes  priority over the others.   For example, one theorem  says that two triangles $\bigtriangleup ABC$ and $ \bigtriangleup A'B'C'  $ are similar if  the length of their corresponding side are proportional. This means there exists a positive number $k>0$ such that 
\[ k=\frac{\overline{AB}}{\overline{A'B'}}=\frac{\overline{AC}}{\overline{A'C'}}=\frac{\overline{BC}}{\overline{B'C'}}. \] 
In such a case,  the   number $k$ is usually called the \textit{ratio of  similarity} (see \cref{Figure1}). 

Note that proportionality of corresponding sides in similar triangles  implies that the corresponding angles are   equal   and vice versa (see \cref{Figure1}). So,  some authors  use the equality of angles to define the similarity of triangles:   two triangles $\bigtriangleup ABC$ and $ \bigtriangleup A'B'C'  $ are similar if
\[ \angle A = \angle A',~  \angle B = \angle B',~ \angle C = \angle C'.\]

Another definition for similarity of triangles says that   two    triangles $\bigtriangleup ABC$ and $ \bigtriangleup A'B'C'  $ are similar if two sides are proportional and the angles between these two sides are equal. It can be shown that the above-mentioned definitions of similarity are equivalent (see \cite{Sib98}).

\section{Intercept Theorem} 
Another concept, which has a close relation with the similarity of triangles, is the   \textit{intercept theorem}.  This theorem  is usually attributed to the Greek philosopher  Thales of Miletus  (circa 624--545 BC)   and sometimes called the \textit{Thales's intercept theorem}. Let us explain this elementary theorem in the following example.

\begin{figure}[H]
	\centering
	\includegraphics[scale=1]{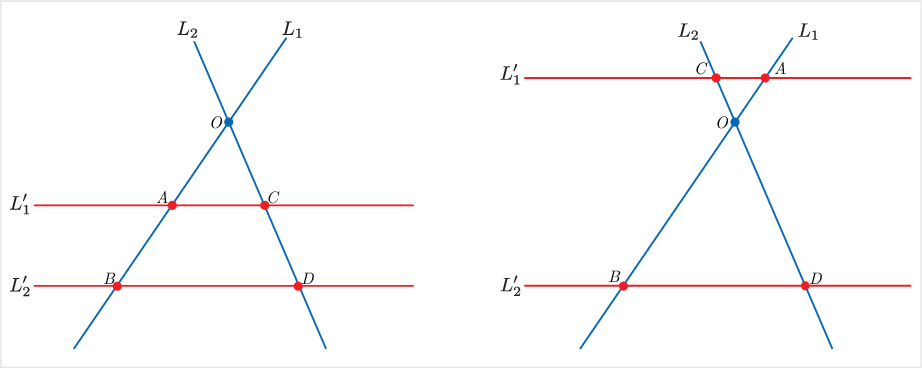}
	%\vspace{-1mm}
	\caption{Intercept theorem}
	\label{Figure2}
\end{figure}

 Consider two straight lines $L_1$ and $L_2$ intersecting in a point $O$ and assume that two parallel lines $L'_1$ and $L'_2$ intersect  $L_1$ and $L_2$ such that they do not pass through $O$. There are only two cases  (see \cref{Figure2}):  \\
 (1) $O$ is in the region bounded by $L'_1$ and $L'_2$, or\\
 (2)  $O$ is not in the region bounded by $L'_1$ and $L'_2$.\\
  If $L_1\cap L'_1=\{A\} $, $L_1\cap L'_2=\{B\} $, $L_2\cap L'_1=\{C\} $, 
and $L_2\cap L'_2=\{D\} $, then the theorem says that
\[ \frac{\overline{OA}}{\overline{OB}}= \frac{\overline{OC}}{\overline{OD}} = \frac{\overline{AC}}{\overline{BD}}. \]

Note that in \cref{Figure2}  the  two triangles $\bigtriangleup AOC $ and $\bigtriangleup BOD $ formed by the  four lines $L_1, L_2, L'_1, L'_2$ are similar.  In fact, it can be shown that this theorem is somehow equivalent to    the similarity of triangles. In other words, by assuming one assertion, one can prove the other. The first proof of this theorem seems to be provided by Euclid in his famous book \textit{Elements} (see \cite{Hea56}). 

\section{Transversals}
\begin{figure}[H]
	\centering
	\includegraphics[scale=1]{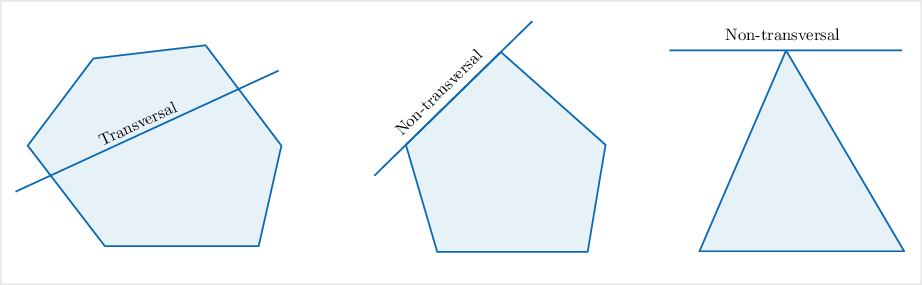}
	%\vspace{-1mm}
	\caption{Transversal and non-transversal lines}
	\label{SMT18-5-e}
\end{figure}

In   elementary geometry, a \textit{transversal line} or a \textit{transversal} is a line intersecting two different lines in the plane in two distinct points. For example, in \cref{Figure2} line $L'_1$ is a transversal with respect to   intersecting lines $L_1$ and $L_2$.  For two dimensional figures such as polygons, a transversal is usually   a line dividing  the figure into two parts whose areas are positive. So it has to intersect at least two sides of the figure. \cref{SMT18-5-e} shows three lines but one of which is   transversal. 

\begin{figure}[H]
	\centering
	\includegraphics[scale=1]{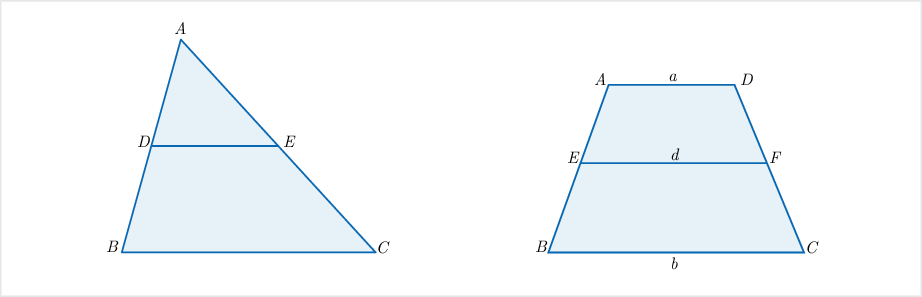}
	%\vspace{-1mm}
	\caption{Transversals of polygons}
	\label{Figure5}
\end{figure}

Of all  possible transversals for a polygon,   the  one   parallel to a third side of the polygon is of special interest.  Some authors consider these special intersecting lines as transversals of polygons. For example, a transversal of a triangle is a line intersecting its two sides and parallel to the third one  (see \cref{Figure5}, left).   In such cases, one can use the intercept theorem to compute the length of a  part of the triangle with respect to the other parts. For a trapezoid, a transversal is usually a line parallel to the two bases which intersects the two legs (\cref{Figure5}, right).

\section{Applications of Intercept Theorem}
Although  intercept theorem seems   elementary, it has many practical applications. This theorem has been   used  for a   long time by mathematicians and surveyors to measure   distances   that were not  capable of measurement using standard methods. 

\begin{figure}[H]
	\centering
	\includegraphics[scale=1]{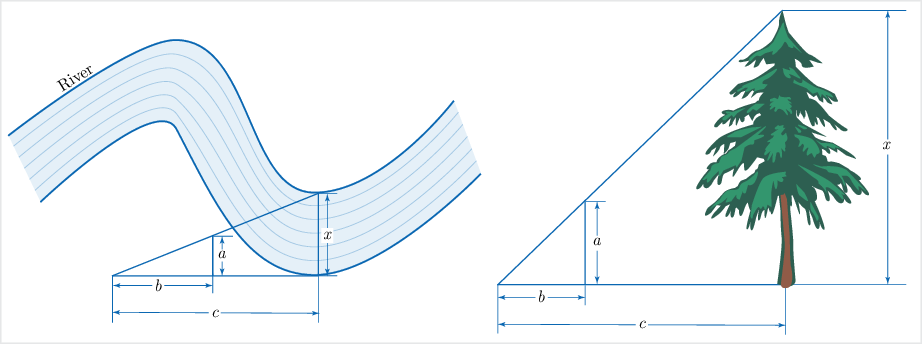}
	%\vspace{-1mm}
	\caption{Applications of intercept theorem}
	\label{SMT18-7-e}
\end{figure}

For example, the theorem can be used to compute  the width of rivers or the height of tall trees or structures. These two situations are depicted in \cref{SMT18-7-e}.   In the case of the tree,   $c$ can be chosen as the length of shadow of the tree and $a$ and $b$ can be the lengths  of  a stick  and its shadow respectively. In both cases, one can easily determine the   values of  reachable quantities $a,b,c$ and then use the intercept theorem to compute the unreachable value $x$ by $x=\frac{ac}{b} $.  

It has been said that  Thales applied the intercept theorem in a similar manner  to measure the height of the \textit{pyramid of Cheops}\footnote{The  pyramid of Cheops (also known as the pyramid of Khufu or the great pyramid of Giza) is the oldest and largest of the pyramids in the Giza pyramid complex bordering present-day Giza in the Greater Cairo, Egypt.} (see \cite{Hea21}).

 \section{Similarity in the SMT}
 As we said before, the   similarity of triangles and transversals were used in some texts of the \textbf{SMT}. We have discussed elsewhere the transversal bisectors of trapezoids, which are the subjects of \textbf{SMT No.\,23} and \textbf{SMT No.\,25} (see \cite{HM22-2}). The main idea of the problems in those texts was to use a transversal parallel to the bases in order to divide the trapezoid into two subtrapezoids with equal areas (see \cref{Figure5}). 
 
 The key property of the transversal of a trapezoid is that its length depends only on the length of the two bases. In fact, if $a,b$ are the lengths of the two bases, then the  length $d$ of the transversal is obtained by
 \begin{equation}\label{eq-a}
 	d = \sqrt{\dfrac{a^2+b^2}{2}}.
 \end{equation}
 The key point to prove formula \cref{eq-a} is to consider the height of the trapezoid and  use the similarity of triangles (see \cite{HM22-2}, Section 3).

\subsection{\textbf{SMT No.\,18}}

\subsubsection*{Transliteration}

\begin{note1}
	\underline{Obverse:  Lines 1-12}\\
	(L1)\hspace{2mm} u\v{s} ki-ta \textit{a-na} u\v{s} an-ta nigin-\textit{ma} $<$10$>$\\
	(L2)\hspace{2mm} a-\v{s}\`{a} an-ta \textit{a-na} a-\v{s}\`{a} [ki]-ta nigin 36\\
	(L3)\hspace{2mm} [sag an]-ta nigin dal nigin ul-gal 20,24\\
	(L4)\hspace{2mm} za-e 36 \textit{\v{s}\`{a}} a-\v{s}\`{a} ki a-\v{s}\`{a} nigin\\
	(L5)\hspace{2mm} \textit{a-na} 4 \textit{a-li-ik-ma} 2,24 \textit{ta-mar}
	\\
	(L6)\hspace{2mm} igi-10 \textit{\v{s}\`{a} mu} ki \textit{mu} nigin \textit{pu-\c{t}\'{u}-}$<$\textit{\'{u}r}$>$ 6 \textit{ta-mar}\\
	(L7)\hspace{2mm} 2,24 \textit{a-na} 6 \textit{i-\v{s}\'{i}-ma} 14,24 \textit{ta-mar}\\
	(L8)\hspace{2mm} 14,24 nigin 3,27,21,[36] \textit{ta-mar a-na} [2] \textit{i-\v{s}\'{i}} [6,54],43,12 \textit{ta-mar}\\
	(L9)\hspace{2mm} [\textit{tu-\'{u}}]\textit{r} 14,24 \textit{a-na} 2 \textit{i-\v{s}\'{i}} [28,48 \textit{ta-mar} $\cdots$] $\cdots$\\

	\underline{Reverse:  Lines 1-3}\\
	(L1)\hspace{2mm} $\cdots$ [$\cdots$ $\cdots$ $\cdots$]\\
	(L2)\hspace{2mm} 30 \textit{ta}-[\textit{mar} $\cdots$ $\cdots$]
	\\
	(L3)\hspace{2mm} \textit{i-\v{s}\'{i}-ma} 20 \textit{ta-ma}[\textit{r} $\cdots$ $\cdots$]
	
\end{note1}

\subsubsection*{Translation}

\underline{Obverse:  Lines 1-9}
\begin{tabbing}
	\hspace{14mm} \= \kill 
	(L1)\> \tabfill{I multiplied the lower length by the upper length, and (the result is) 10,0.}\index{length}\\
	(L2)\> \tabfill{I multiplied the upper area by the lower area, (and the result is) 36,0,0.}\index{area}\\
	(L3)\> \tabfill{I added the squared upper width (and) the squared transversal\index{transversal line}, (and the result is) 20,24.}\index{width}\\
	(L4)\> \tabfill{You, 36,0,0 that is the (result of) multiplication of the area by (another) area,}\index{area}\\
	(L5)\> \tabfill{multiply (it) by 4, and you see 2,24,0,0.}\\
	(L6)\> \tabfill{Make the reciprocal of 10,0 that is (the result of) multiplication of the perpendicular (that is, the lower length) by the perpendicular (that is, the upper length), (and) you see 0;0,6.}\index{reciprocal of a number}\index{perpendicular line}\index{length} \\ 
	(L7)\> \tabfill{Multiply 2,24,0,0 by 0;0,6, and you see 14,24.}\\ 
	(L8)\> \tabfill{Square 14,24, (and) you see 3,27,21,36. Multiply (it) by 2, (and) you see 6,54,43,12.}\\
	(L9)\> \tabfill{Return. Multiply 14,24 by 2, (and) you see 28,48. $\cdots$ $\cdots$.}
\end{tabbing}

\noindent
\underline{Reverse:  Lines 1-3}
\begin{tabbing}
	\hspace{14mm} \= \kill 
	(L1)\> \tabfill{$\cdots$ $\cdots$ $\cdots$.}\\
	(L2)\> \tabfill{you see 30. $\cdots$ $\cdots$.}\\
	(L3)\> \tabfill{Multiply (30) by (0;40), and you see 20.} 
	
\end{tabbing}

\subsubsection*{Mathematical Interpretation}
The Susa scribe   is considering the dimensions and areas of a right triangle  and a trapezoid  as   shown in the following figure. We have a right triangle with a transversal line which is  parallel to the  base of the right triangle  and  divides it into two figures: a smaller right triangle  and a right trapezoid.

\begin{figure}[H]
	\centering
	\includegraphics[scale=1]{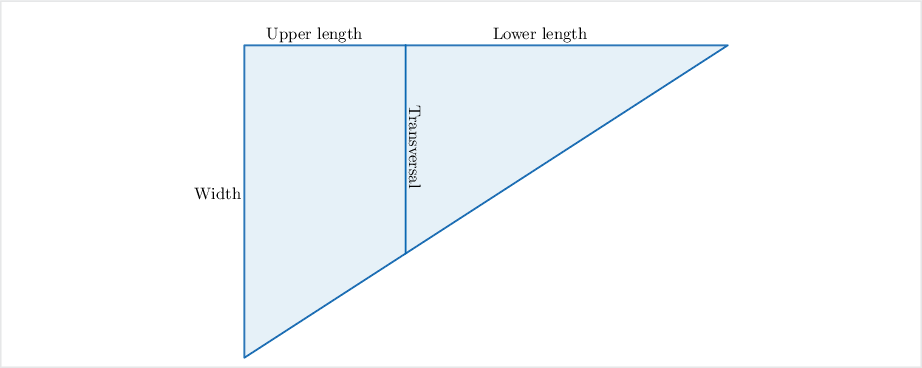}
	%\vspace{-1mm}
	\caption{Transversal of a right triangle}
	\label{Figure7}
\end{figure}

According to  \cref{Figure7} and the translation, there are four variables in this problem: the upper length, the lower length, the width  and the transversal.  We use the following symbols for these variables in our  discussion\footnote{Following Babylonian tradition, for determining the lower and the upper lengths, we consider right-to-left direction, while for the lower and upper width we take the down-to-up direction.}: 

\[ \begin{cases}
	x=\text{the upper length},\\
	y=\text{the lower length},\\
	z=\text{the width\index{width}},\\
	w=\text{the transversal}.
\end{cases} \]\index{transversal line}\index{length}

\noindent
(Note that the scribe has implicitly assumed that $z>w$.)  

\begin{figure}[H]
	\centering
	\includegraphics[scale=1]{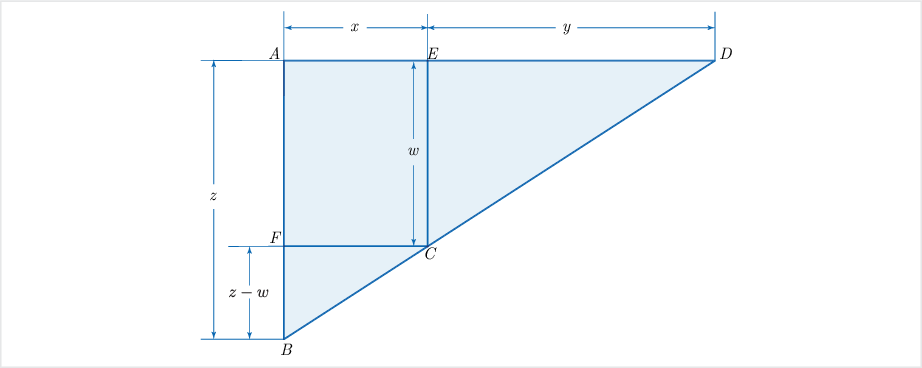}
	%\vspace{-1mm}
	\caption{Dimensions of a right triangle with transversal}
	\label{Figure8}
\end{figure}

In \cref{Figure8}, we  consider the perpendicular lines  from $C$ onto $AB$ and  $AD$ and denote the intersection points by $F$ and $E$. As shown in   \cref{Figure8}, we have labeled the vertices and    parts of the figure as follows:
\[ \overline{AE}=x,~\overline{ED}=y,~\overline{AB}=z,~\overline{EC}=w. \]
In this case,   we have $\overline{FB}=z-w$.

If we use this notation, lines 1-3 of the text  provide us with the following system of simultaneous equations:

\begin{equation}\label{equ-SMT18-a}
	\begin{dcases}
		xy=10,0\\
		\left(\dfrac{1}{2}x(z+w)\right)\times \left(\dfrac{1}{2}yw\right)=36,0,0\\
		z^2+w^2=20,24.
	\end{dcases}
\end{equation} 
Note that the expression $\frac{1}{2}x(z+w) $ is the very area of the right trapezoid  $ABCE$ and the other expression $ \frac{1}{2}yw $ is that of the right triangle  $\triangle CDE$.

As we can see, there are four variables and three equations which makes it impossible to solve the system \cref{equ-SMT18-a}. Because of this difficulty, the scribe   solved this system of equations in two steps. In the first  step, he eliminates $x$ and $y$ to  obtain a system of equations with respect to only  $w$ and $z$ in order to find their values. Then, in the second step, he tries to find  the values of $x$ and $y$ by using the first equation in \cref{equ-SMT18-a} and a property of the figure which provides him with one more equation  with respect to $x$ and $y$ (this is where he is   using the   intercept theorem). We describe each step in details.  \\

\noindent
\textbf{Step 1.} 
According to lines 4-7, the scribe  uses the first and the second equations in \cref{equ-SMT18-a} to  get the value of $w(z+w)$ as follows:

\begin{align*}
	&~~  \left(\dfrac{1}{2}x(z+w)\right)\times \left(\dfrac{1}{2}yw\right)=36,0,0 \\
	\Longrightarrow~~&~~   4\times \left(\dfrac{1}{2}x(z+w)\right)\times \left(\dfrac{1}{2}yw\right)=4\times(36,0,0)  \\
	\Longrightarrow~~&~~    xy \times  w(z+w) =  2,24,0,0  \\
	\Longrightarrow~~&~~    (10,0)\times  w(z+w) =  2,24,0,0  \\
	\Longrightarrow~~&~~    w(z+w)=  \frac{1}{(10,0)} \times(2,24,0,0) \\
	\Longrightarrow~~&~~    w(z+w)= (0;0,6)\times(2,24,0,0)
\end{align*} 
thus
\begin{equation}\label{equ-SMT18-b}
	w(z+w) =14,24.
\end{equation} 
According to line 8, by squaring and then doubling both sides of   \cref{equ-SMT18-b}, we get 
\begin{align*}
	&~~  w(z+w) =14,24 \\
	\Longrightarrow~~&~~   \Big(w(z+w)\Big)^2 =(14,24)^2  \\
	\Longrightarrow~~&~~   w^2(z+w)^2 =3,27,21,36 \\
	\Longrightarrow~~&~~    2w^2(z+w)^2  =2\times(3,27,21,36)  \\
	\Longrightarrow~~&~~    2w^2(z+w)^2 =6,54,43,12.
\end{align*}
So
\begin{equation}\label{equ-SMT18-c}
	 2w^2  (z+w)^2  =6,54,43,12.
\end{equation}
Next, according to line 9, we multiply both sides of \cref{equ-SMT18-b} by 2 to obtain
\begin{equation}\label{equ-SMT18-ca}
	2w(z+w) =28,48.
\end{equation} 

It seems that at this point of the text  the scribe has used new variables to proceed. We may recover the next calculations as follows. First introduce new variables $X$ and $Y$ by
\begin{equation}\label{equ-SMT18-d}
	\begin{dcases}
		X=(z+w)^2\\
		Y=2w^2.
	\end{dcases}
\end{equation} 
It follows from \cref{equ-SMT18-c} and \cref{equ-SMT18-d} that
\begin{equation}\label{equ-SMT18-e}
	XY=6,54,43,12.
\end{equation} 
On the other hand, from \cref{equ-SMT18-a}, \cref{equ-SMT18-ca} and \cref{equ-SMT18-d} we can write
\begin{align*}
	X+Y&=(z+w)^2+2w^2 \\
	&=z^2+w^2+2zw+2w^2\\
	&=z^2+w^2+2w(z+w)\\
	&=20,24+28,48\\
	&=49,12
\end{align*}
thus

\begin{equation}\label{equ-SMT18-f}
	X+Y=49,12.
\end{equation} 
Now, we can apply the usual Babylonian method, i.e., completing the square    to find the values of $X$ and $Y$ satisfying simultaneous equations  \cref{equ-SMT18-e} and \cref{equ-SMT18-f}. Since
\[ \dfrac{X+Y}{2}=\frac{49,12}{2}=24,36 \]
we  can write
\begin{align*}
	\frac{X-Y}{2}&=\sqrt{\left(\frac{X+2}{2}\right)^2-XY}\\
	&=\sqrt{\left(\frac{49,12}{2}\right)^2-6,54,43,12}\\
	&=\sqrt{(24,36)^2-6,54,43,12}\\
	&=\sqrt{10,5,9,36-6,54,43,12}\\
	&=\sqrt{3,10,26,24}\\
	&=\sqrt{(13,48)^2}\\
	&=13,48
\end{align*}
thus we obtain\footnote{Note that the scribe  might have   computed   $ \sqrt{3,10,26,24}=\sqrt{2^4\times 3^4\times 23^2}=2^2\times 3^2\times 23^1=13,48$.}
\begin{equation}\label{equ-SMT18-g}
	\frac{X-Y}{2}=13,48.
\end{equation} 
It follows from the \cref{equ-SMT18-f} and \cref{equ-SMT18-g} that
\[ X=\frac{X+Y}{2}+\frac{X-Y}{2}=\frac{49,12}{2}+13,48=24,36+13,48=38,24 \]
and
\[ Y=\frac{X+Y}{2}-\frac{X-Y}{2}=\frac{49,12}{2}-13,48=24,36-13,48=10,48. \]
Therefore, we obtain
\begin{equation}\label{equ-SMT18-h}
	X=38,24\ \ \ \text{and} \ \ \ Y=10,48.
\end{equation} 
Now, we can use \cref{equ-SMT18-d} and  \cref{equ-SMT18-h} to compute the values of $z$ and $w$ as follows:
\begin{align*}
	&~~   2w^2 =10,48 \\
	\Longrightarrow~~&~~      w^2   = \frac{10,48}{2}  \\
	\Longrightarrow~~&~~    w^2   = 5,24  \\
	\Longrightarrow~~&~~    w = \sqrt{5,24}  \\
	\Longrightarrow~~&~~    w = \sqrt{(18)^2}  \\
	\Longrightarrow~~&~~    w=18 
\end{align*}
and
\begin{align*}
	&~~   (z+w)^2 =38,24 \\
	\Longrightarrow~~&~~     z+w   =\sqrt{38,24}  \\
	\Longrightarrow~~&~~     z+w   =\sqrt{(48)^2}  \\
	\Longrightarrow~~&~~    z+18  =48 \\
	\Longrightarrow~~&~~    z=48-18  \\
	\Longrightarrow~~&~~    z=30.
\end{align*}
Therefore, we get
\begin{equation}\label{equ-SMT18-i}
	w=18\ \ \ \text{and}\ \ \ z=30
\end{equation} 
which completes the first step.\\

\noindent
\textbf{Step 2.} 
In the second step, the scribe needs to use another condition   in order to find an equation involving only $x$ and $y$ other than  the first equation $xy=10,0$. In fact, if we substitute the values of $w=18$ and $z=30$ into the second equation of \cref{equ-SMT18-a} and simplify,  it is clear  that we get the first equation $xy=10,0$.

A second equation  can be obtained from the properties of the transversal line  in   \cref{Figure8}.  Since $AB\parallel EC $ and $AD \parallel CF$, the   intercept theorem implies that two triangles $  \triangle CDE $ and $  \triangle BCF $ are similar and thus 
\[ \dfrac{\overline{CF}}{\overline{BF}}=\dfrac{\overline{DE}}{\overline{CE}} \] 
or equivalently
\begin{equation}\label{equ-SMT18-j}
	\frac{x}{z-w} =\frac{y}{w}.
\end{equation} 
It follows from  \cref{equ-SMT18-i} and \cref{equ-SMT18-j} that
\[ \frac{x}{12} =\frac{y}{18} \]
or  
\begin{equation}\label{equ-SMT18-k}
	x=\frac{2}{3}y.
\end{equation} 

Equation  \cref{equ-SMT18-k} is   the very condition    that the scribe  has used  to finish the solution. Thus, we have obtained the following system of equations with respect to $x$ and $y$ only:
\begin{equation}\label{equ-SMT18-ka}
	\begin{dcases}
		xy=10,0\\
		x=\frac{2}{3}y.
	\end{dcases}	
\end{equation}

Let us solve this system of equations. By substituting the value of $x$ with respect to $y$ given   by the second  equation of \cref{equ-SMT18-ka} into the first equation,  we can write
\begin{align*}
	&~~   xy=10,0 \\
	\Longrightarrow~~&~~     \left(\frac{2}{3}y\right)y=10,0  \\
	\Longrightarrow~~&~~    \frac{2}{3}y^2=10,0  \\
	\Longrightarrow~~&~~    y^2=\frac{3}{2}\times (10,0) \\
	\Longrightarrow~~&~~    y^2= 15,0  \\
	\Longrightarrow~~&~~    y = \sqrt{15,0}\\
	\Longrightarrow~~&~~    y = 30.   
\end{align*}
So, according to line 2 on the reverse, we get  $y=30$. Finally,  according to line 3 on the reverse,  we have
\[ x=\frac{2}{3}y=\frac{2}{3}\times 30=20.  \]
Thus  the solutions of the system of equations  \cref{equ-SMT18-ka}   are $x=20$ and $y=30$. Ultimately, the solutions of the main system of equation \cref{equ-SMT18-a} 
are given by
\[x=20,\ \ \ y=30, \ \ \ z=30,\ \ \ w=18. \]

\begin{remark}\label{rem-SMT18-a}
	A  mathematical interpretation of this text has been given by Friberg\index{J\"{o}ran Friberg} in \cite{Fri07-2}.
\end{remark}

\subsection{Conclusion}\label{SS-C-SMT18} 
From our mathematical interpretation  of \textbf{SMT No.\,18} we consider it is apparent that the Susa scribes were familiar with the idea of two similar triangles and the relation between their sides. This observation is of great importance to the    history of mathematics  in that it confirms    the origin of similarity and intercept theorem date to approximately a millennium    before the Greeks. 
 
\textbf{SMT No.\,18} also deals with one of the most complicated systems of equations in Babylonian mathematics because there are four unknown variables  involved in the system. Although there are only three equations  and four unknowns at the beginning of the problem, the Susa scribe     has used his geometrical knowledge to obtain another equation employing similarity of triangles in order to achieve the solution.


{\small 
\begin{thebibliography}{00000000}

\bibitem[BR61]{BR61} 
E. M. Bruins and M. Rutten,    \textit{Textes Math\a'{e}matiques de Suse} [\textbf{TMS}], Librairie Orientaliste Paul Geuthner,   Paris, 1961.


\bibitem[Fri07-2]{Fri07-2}
{J. Friberg}, {\em Amazing Traces of a Babylonian Origin in Greek Mathematics}, World Scientific, 2007.

 


\bibitem[Hea21]{Hea21}
{S. T. Heath}, {\em A History of Greek Mathematics: Volumes I-II}, Oxford University Press,   1921.


\bibitem[Hea56]{Hea56}
{S. T. Heath}, {\em The Thirteen Books of the Elements, Vol. 2}, Dover Publications, Second Edition,   1956.

\bibitem[HM22-2]{HM22-2}
{N. Heydari and K. Muroi}, {\em Bisection of Trapezoids in Elamite Mathematics}, arXiv:2211.11139 [math.HO], November 2022, 21 pages: \url{https://arxiv.org/abs/2211.11139}.




	\bibitem[Sib98]{Sib98}
{T. Q. Sibley}, {\em $\pi$: The Geometric Viewpoint : a Survey of Geometries}, Addison-Wesley,  1998.




 \end{thebibliography} 
 
 }
\end{document}